\documentclass[12pt]{article}

\usepackage{amsmath,amsthm}
\usepackage{amssymb}
\usepackage{color}
\usepackage{xspace}
\usepackage{pstricks,pst-node}
\usepackage[colorlinks=true,
linkcolor=green,
filecolor=brown,
citecolor=green]{hyperref}
\definecolor{gray1}{rgb}{0.7,0.7,0.7}
\definecolor{black}{rgb}{0,0,0}
\definecolor{maggy}{rgb}{0.6,0.2,1}
\definecolor{browne}{rgb}{0.8,0,.2}

\def\red{\textcolor{red} }
\def\blue{\textcolor{blue} }
\def\gray1{\textcolor{gray1} }
\def\black{\textcolor{black} }
\def\maggy{\textcolor{maggy} }
\def\browne{\textcolor{browne} }

\newrgbcolor{brown}{0.8 0 0.2}

\def\IOTs{increasing ordered trees\xspace}
\def\IOT{increasing ordered tree\xspace}
\def\v{\vert}
\def\s{\ensuremath{\mathcal S}\xspace}
\def\a{\ensuremath{\mathcal A}\xspace}
\def\d{\ensuremath{\mathcal D}\xspace}

\def\u{\ensuremath{\mathcal U}\xspace}

\def\I{\ensuremath{\mathcal I}\xspace} 
\def\J{\ensuremath{\mathcal J}\xspace}  
\def\p{\ensuremath{\mathcal P}\xspace}

\def\w{\ensuremath{\mathcal W}\xspace}
\def\x{\ensuremath{\mathcal X}\xspace}
\def\y{\ensuremath{\mathcal Y}\xspace}
\def\z{\ensuremath{\mathcal Z}\xspace}

\def\si{\sigma}

\def\mbf#1{\mathchoice{\hbox{\boldmath $\displaystyle #1$}}
        {\hbox{\boldmath $\textstyle #1$}}
        {\hbox{\boldmath $\scriptstyle #1$}}
        {\hbox{\boldmath $\scriptscriptstyle #1$}}} 

\catcode`\@=11

\thicklines
\newskip\Einheit \Einheit=0.6cm
\newcount\xcoord \newcount\ycoord
\newdimen\xdim \newdimen\ydim \newdimen\PfadD@cke \newdimen\Pfadd@cke
\def\PfadDicke#1{\PfadD@cke#1 \divide\PfadD@cke by2 
\Pfadd@cke\PfadD@cke \multiply\PfadD@cke by2}
\long\def\LOOP#1\REPEAT{\def\BODY{#1}\ITERATE}
\def\ITERATE{\BODY \let\next\ITERATE \else\let\next\relax\fi \next}
\let\REPEAT=\fi
\def\Punkt{\hbox{\raise-2pt\hbox to0pt{\hss\scriptsize$\bullet$\hss}}}

\def\DuennPunkt(#1,#2){\unskip
  \raise#2 \Einheit\hbox to0pt{\hskip#1 \Einheit
          \raise-1.5pt\hbox to0pt{\hss\tiny$\bullet$\hss}\hss}}
		  
\def\DotPunkt(#1,#2){\unskip
  \raise#2 \Einheit\hbox to0pt{\hskip#1 \Einheit
          \raise-1.5pt\hbox to0pt{\hss\normalsize$\dot$\hss}\hss}}

\def\NormalPunkt(#1,#2){\unskip
  \raise#2 \Einheit\hbox to0pt{\hskip#1 \Einheit
          \raise-3pt\hbox to0pt{\hss\large$\bullet$\hss}\hss}}
\def\DickPunkt(#1,#2){\unskip
  \raise#2 \Einheit\hbox to0pt{\hskip#1 \Einheit
          \raise-4pt\hbox to0pt{\hss\Large$\bullet$\hss}\hss}}
\def\Kreis(#1,#2){\unskip
  \raise#2 \Einheit\hbox to0pt{\hskip#1 \Einheit
          \raise-4pt\hbox to0pt{\hss\Large$\circ$\hss}\hss}}
\def\Diagonale(#1,#2)#3{\unskip\leavevmode
  \xcoord#1\relax \ycoord#2\relax
      \raise\ycoord \Einheit\hbox to0pt{\hskip\xcoord \Einheit
         \unitlength\Einheit
         \line(1,1){#3}\hss}}
\def\AntiDiagonale(#1,#2)#3{\unskip\leavevmode
  \xcoord#1\relax \ycoord#2\relax \advance\xcoord by -0.05\relax
      \raise\ycoord \Einheit\hbox to0pt{\hskip\xcoord \Einheit
         \unitlength\Einheit
         \line(1,-1){#3}\hss}}
\def\Pfad(#1,#2),#3\endPfad{\unskip\leavevmode
  \xcoord#1 \ycoord#2 \thicklines\ZeichnePfad#3\endPfad\thinlines}
\def\ZeichnePfad#1{\ifx#1\endPfad\let\next\relax
  \else\let\next\ZeichnePfad
    \ifnum#1=1
      \raise\ycoord \Einheit\hbox to0pt{\hskip\xcoord \Einheit
         \vrule height\Pfadd@cke width1 \Einheit depth\Pfadd@cke\hss}%
      \advance\xcoord by 1
    \else\ifnum#1=2
      \raise\ycoord \Einheit\hbox to0pt{\hskip\xcoord \Einheit
         \unitlength\Einheit
         \line(0,1){1}\hss}
      \advance\xcoord by 0
      \advance\ycoord by 1
    \else\ifnum#1=3
      \raise\ycoord \Einheit\hbox to0pt{\hskip\xcoord \Einheit
         \unitlength\Einheit
         \line(1,1){1}\hss}
      \advance\xcoord by 1
      \advance\ycoord by 1
    \else\ifnum#1=4
      \raise\ycoord \Einheit\hbox to0pt{\hskip\xcoord \Einheit
         \unitlength\Einheit
         \line(1,-1){1}\hss}
      \advance\xcoord by 1
      \advance\ycoord by -1
	  \else\ifnum#1=5
      \raise\ycoord \Einheit\hbox to0pt{\hskip\xcoord \Einheit
         \unitlength\Einheit
         \line(2,1){2}\hss}
      \advance\xcoord by 2
      \advance\ycoord by 1
	  \else\ifnum#1=6
      \raise\ycoord \Einheit\hbox to0pt{\hskip\xcoord \Einheit
         \unitlength\Einheit
         \line(2,-1){2}\hss}
      \advance\xcoord by 2
      \advance\ycoord by -1
	  \else\ifnum#1=7
      \raise\ycoord \Einheit\hbox to0pt{\hskip\xcoord \Einheit
         \unitlength\Einheit
         \line(3,1){3}\hss}
      \advance\xcoord by 3
      \advance\ycoord by 1
	  \else\ifnum#1=8
      \raise\ycoord \Einheit\hbox to0pt{\hskip\xcoord \Einheit
         \unitlength\Einheit
         \line(3,-1){3}\hss}
      \advance\xcoord by 3
      \advance\ycoord by -1
    \fi\fi\fi\fi\fi\fi\fi\fi
  \fi\next}
\def\hSSchritt{\leavevmode\raise-.4pt\hbox 
to0pt{\hss.\hss}\hskip.2\Einheit
  \raise-.4pt\hbox to0pt{\hss.\hss}\hskip.2\Einheit
  \raise-.4pt\hbox to0pt{\hss.\hss}\hskip.2\Einheit
  \raise-.4pt\hbox to0pt{\hss.\hss}\hskip.2\Einheit
  \raise-.4pt\hbox to0pt{\hss.\hss}\hskip.2\Einheit}
\def\vSSchritt{\vbox{\baselineskip.2\Einheit\lineskiplimit0pt
\hbox{.}\hbox{.}\hbox{.}\hbox{.}\hbox{.}}}
\def\DSSchritt{\leavevmode\raise-.4pt\hbox to0pt{%
  \hbox to0pt{\hss.\hss}\hskip.2\Einheit
  \raise.2\Einheit\hbox to0pt{\hss.\hss}\hskip.2\Einheit
  \raise.4\Einheit\hbox to0pt{\hss.\hss}\hskip.2\Einheit
  \raise.6\Einheit\hbox to0pt{\hss.\hss}\hskip.2\Einheit
  \raise.8\Einheit\hbox to0pt{\hss.\hss}\hss}}
\def\dSSchritt{\leavevmode\raise-.4pt\hbox to0pt{%
  \hbox to0pt{\hss.\hss}\hskip.2\Einheit
  \raise-.2\Einheit\hbox to0pt{\hss.\hss}\hskip.2\Einheit
  \raise-.4\Einheit\hbox to0pt{\hss.\hss}\hskip.2\Einheit
  \raise-.6\Einheit\hbox to0pt{\hss.\hss}\hskip.2\Einheit
  \raise-.8\Einheit\hbox to0pt{\hss.\hss}\hss}}
\def\SPfad(#1,#2),#3\endSPfad{\unskip\leavevmode
  \xcoord#1 \ycoord#2 \ZeichneSPfad#3\endSPfad}
\def\ZeichneSPfad#1{\ifx#1\endSPfad\let\next\relax
  \else\let\next\ZeichneSPfad
    \ifnum#1=1
      \raise\ycoord \Einheit\hbox to0pt{\hskip\xcoord \Einheit
         \hSSchritt\hss}%
      \advance\xcoord by 1
    \else\ifnum#1=2
      \raise\ycoord \Einheit\hbox to0pt{\hskip\xcoord \Einheit
        \hbox{\hskip-2pt \vSSchritt}\hss}%
      \advance\ycoord by 1
    \else\ifnum#1=3
      \raise\ycoord \Einheit\hbox to0pt{\hskip\xcoord \Einheit
         \DSSchritt\hss}
      \advance\xcoord by 1
      \advance\ycoord by 1
    \else\ifnum#1=4
      \raise\ycoord \Einheit\hbox to0pt{\hskip\xcoord \Einheit
         \dSSchritt\hss}
      \advance\xcoord by 1
      \advance\ycoord by -1
    \fi\fi\fi\fi
  \fi\next}
\def\Koordinatenachsen(#1,#2){\unskip
 \hbox to0pt{\hskip-.5pt\vrule height#2 \Einheit width.5pt depth1 
\Einheit}%
 \hbox to0pt{\hskip-1 \Einheit \xcoord#1 \advance\xcoord by1
    \vrule height0.25pt width\xcoord \Einheit depth0.25pt\hss}}
\def\Koordinatenachsen(#1,#2)(#3,#4){\unskip
 \hbox to0pt{\hskip-.5pt \ycoord-#4 \advance\ycoord by1
    \vrule height#2 \Einheit width.5pt depth\ycoord \Einheit}%
 \hbox to0pt{\hskip-1 \Einheit \hskip#3\Einheit 
    \xcoord#1 \advance\xcoord by1 \advance\xcoord by-#3 
    \vrule height0.25pt width\xcoord \Einheit depth0.25pt\hss}}
\def\Gitter(#1,#2){\unskip \xcoord0 \ycoord0 \leavevmode
  \LOOP\ifnum\ycoord<#2
    \loop\ifnum\xcoord<#1
      \raise\ycoord \Einheit\hbox to0pt{\hskip\xcoord 
\Einheit\Punkt\hss}%
      \advance\xcoord by1
    \repeat
    \xcoord0
    \advance\ycoord by1
  \REPEAT}
\def\Gitter(#1,#2)(#3,#4){\unskip \xcoord#3 \ycoord#4 \leavevmode
  \LOOP\ifnum\ycoord<#2
    \loop\ifnum\xcoord<#1
      \raise\ycoord \Einheit\hbox to0pt{\hskip\xcoord 
\Einheit\Punkt\hss}%
      \advance\xcoord by1
    \repeat
    \xcoord#3
    \advance\ycoord by1
  \REPEAT}
\def\Label#1#2(#3,#4){\unskip \xdim#3 \Einheit \ydim#4 \Einheit
  \def\lo{\advance\xdim by-.5 \Einheit \advance\ydim by.5 \Einheit}%
  \def\llo{\advance\xdim by-.25cm \advance\ydim by.5 \Einheit}%
  \def\loo{\advance\xdim by-.5 \Einheit \advance\ydim by.25cm}%
  \def\o{\advance\ydim by.25cm}%
  \def\ro{\advance\xdim by.5 \Einheit \advance\ydim by.5 \Einheit}%
  \def\rro{\advance\xdim by.25cm \advance\ydim by.5 \Einheit}%
  \def\roo{\advance\xdim by.5 \Einheit \advance\ydim by.25cm}%
  \def\l{\advance\xdim by-.30cm}%
  \def\r{\advance\xdim by.30cm}%
  \def\lu{\advance\xdim by-.5 \Einheit \advance\ydim by-.6 \Einheit}%
  \def\llu{\advance\xdim by-.25cm \advance\ydim by-.6 \Einheit}%
  \def\luu{\advance\xdim by-.5 \Einheit \advance\ydim by-.30cm}%
  \def\u{\advance\ydim by-.30cm}%
  \def\ru{\advance\xdim by.5 \Einheit \advance\ydim by-.6 \Einheit}%
  \def\rru{\advance\xdim by.25cm \advance\ydim by-.6 \Einheit}%
  \def\ruu{\advance\xdim by.5 \Einheit \advance\ydim by-.30cm}%
  #1\raise\ydim\hbox to0pt{\hskip\xdim
     \vbox to0pt{\vss\hbox to0pt{\hss$#2$\hss}\vss}\hss}%
}
\catcode`\@=12

\begin{document}
\newtheorem{theorem}{Theorem}
\newtheorem{defn}[theorem]{Definition}
\newtheorem{lemma}[theorem]{Lemma}
\newtheorem{prop}[theorem]{Proposition}
\newtheorem{cor}[theorem]{Corollary}
\begin{center}
{\Large
A Bijection to Count (1-23-4)-Avoiding Permutations                            \\ 
}

\vspace{10mm}
DAVID CALLAN  \\
Department of Statistics,  University of Wisconsin-Madison,\\ 
1300 University Ave, Madison, WI \ 53706-1532  \\
{\bf callan@stat.wisc.edu}  \\
\vspace{3mm}

Aug 10, 2010
\end{center}

\begin{abstract}
A permutation is (1-23-4)-avoiding if it contains no four entries, 
increasing left to right, with the middle two adjacent in the 
permutation. Here we give a 2-variable recurrence for the number of 
such permutations, improving on the previously known 4-variable recurrence. 
At the heart of the proof is a bijection from (1-23-4)-avoiding 
permutations to increasing ordered trees whose leaves, taken in 
preorder, are also increasing.
\end{abstract}

\section{Introduction}

\vspace*{-4mm}

\hspace*{6mm} There is a large literature on pattern avoidance in permutations;  
Mikl\'{o}s B\'{o}na's book \cite{bonabook} contains a good bibliography. 
Babson and Steingr\'{\i}msson \cite{steinGrim} introduced the notion 
of ``dashed'' pattern: the absence of a dash indicates the 
corresponding letters in the permutation must be adjacent in the 
permutation. Thus a 
permutation avoids the dashed pattern 1-23-4 if it contains no four entries, 
increasing left to right, with the middle two adjacent. Dashed patterns are also known as vincular patterns and are subsumed by the notion of bivincular pattern \cite{dukes10} where permutation 
entries may be required to be consecutive in value and/or adjacent 
in position. Sergi Elizalde \cite{elizalde06} obtained 
aymptotic bounds for the number $u(n)$ of  permutations of 
$[n]:=\{1,2,\ldots,n\}$ that avoid 1-23-4.
A complicated recurrence for $u(n)$ is posted to the 
\htmladdnormallink{On-Line Encyclopedia of Integer 
Sequences}{http://www.research.att.com:80/~njas/sequences/Seis.html} at
\htmladdnormallink{A113227}{http://www.research.att.com:80/cgi-bin/access.cgi/as/njas/sequences/eisA.cgi?Anum=A113227}.
The purpose of this paper is to establish a much simpler recurrence: 
$u(n)=\sum_{k=1}^{n}u(n,k)$ where $u(n,k)$ is defined by the recurrence
\begin{equation}
\hspace*{15mm} u(n,k)=u(n-1,k-1)+k \sum_{j=k}^{n-1}u(n-1,j) \hspace*{15mm} 1\le k \le n
    \label{recur}
\end{equation}
with initial conditions $u(0,0)=1$ and $u(n,0)=0$ for $n\ge 1$.

In Section 2, we show that this recurrence counts $n$-edge \IOTs 
with increasing leaves, by outdegree $k$ of the root. This reduces the problem to identification of these trees and (1-23-4)-avoiding permutations, our main bijection. In 
Section 3, we review Stirling and Gessel permutations and give a baby version of the main bijection. In 
Section 4, we introduce what we call 2-configurations, 
a generalization of the notion of set partition in which some 
entries are allowed to appear twice.
Sections 5 and 6 present the main bijection, between (1-23-4)-avoiding permutations 
of $[n]$ and $n$-edge \IOTs with increasing leaves, in the process 
representing both the permutations and the trees as 2-configurations. In 
Section 7, we observe that (1-23-4)-avoiding permutations are also 
equinumerous with certain marked-up Dyck paths. Section 8 gives a 4-variable generating function for statistics on permutations that arise naturally earlier in the paper. Finally, Section 9 briefly raises some further questions

\section{Increasing ordered trees with increasing leaves}\label{iot}

\vspace*{-4mm}

\hspace*{6mm} An  \emph{increasing ordered tree} of size $n$  
is an ordered tree (sometimes called 
a plane tree) with $n+1$ labeled vertices, the standard label set being  
$0,1,2,\ldots,n$, such that each child exceeds its parent. Thus  
size measures the number of edges. 
Standard \IOTs of size $n$ are counted by the odd double 
factorial $(2n-1)!!$ (see, e.g.,\,\cite{stirpoly,dblfac}).
An \IOT has \emph{increasing leaves} if its leaves, taken in preorder 
(or ``walkaround'' order) are increasing. The figure below shows two \IOTs, 
the first has increasing leaves while the second does not.

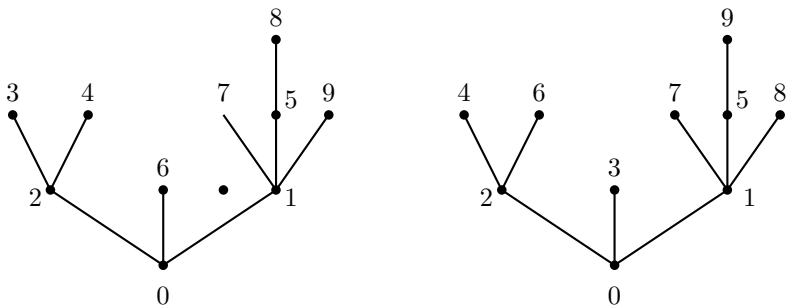
\begin{center} 

\begin{pspicture}(-5,-0.3)(5,3.4)
\psset{unit=1cm}

\psline(-5,2)(-4.5,1)(-4,2)
\psline(-4.5,1)(-3,0)(-3,1)
\psline(-3,0)(-1.5,1)(-1.5,2)(-1.5,3)
\psline(-2.2,2)(-1.5,1)(-0.8,2)

\psdots(-5,2)(-4.5,1)(-4,2)(-3,0)(-3,1)(-1.5,1)(-1.5,2)(-1.5,3)(-2.2,1)(-0.8,2)

\rput(-3,-0.4){\textrm{{\footnotesize $0$}}}
\rput(-4.7,.9){\textrm{{\footnotesize $2$}}}
\rput(-5,2.3){\textrm{{\footnotesize $3$}}}
\rput(-1.3,.9){\textrm{{\footnotesize $1$}}}
\rput(-4,2.3){\textrm{{\footnotesize $4$}}}
\rput(-3,1.3){\textrm{{\footnotesize $6$}}}
\rput(-2.2,2.3){\textrm{{\footnotesize $7$}}}
\rput(-1.3,2.2){\textrm{{\footnotesize $5$}}}
\rput(-1.5,3.3){\textrm{{\footnotesize $8$}}}
\rput(-0.8,2.3){\textrm{{\footnotesize $9$}}}

\psline(1,2)(1.5,1)(2,2)
\psline(1.5,1)(3,0)(4.5,1)(3.8,2)
\psline(4.5,3)(4.5,2)(4.5,1)(5.2,2)
\psline(3,1)(3,0)

\psdots(1,2)(1.5,1)(2,2)(3,0)(3,1)(4.5,1)(3.8,2)(4.5,3)(4.5,2)(5.2,2)

\rput(3,-0.4){\textrm{{\footnotesize $0$}}}
\rput(1.3,.9){\textrm{{\footnotesize $2$}}}
\rput(3,1.3){\textrm{{\footnotesize $3$}}}
\rput(4.8,.9){\textrm{{\footnotesize $1$}}}
\rput(1,2.3){\textrm{{\footnotesize $4$}}}
\rput(2,2.3){\textrm{{\footnotesize $6$}}}
\rput(3.8,2.3){\textrm{{\footnotesize $7$}}}
\rput(4.7,2.2){\textrm{{\footnotesize $5$}}}
\rput(5.2,2.3){\textrm{{\footnotesize $8$}}}
\rput(4.5,3.3){\textrm{{\footnotesize $9$}}}

\end{pspicture}
\end{center}

Now let $\I_{n,k}$ denote the set of size-$n$ \IOTs with 
increasing leaves in which the root has $k$ children and let $u(n,k)
=\v\, \I_{n,k}\,\v$. Consider a tree in $\I_{n,k}$.
If vertex 1 is a leaf, then it must be the leftmost child of the root 
(else the leaves would not be increasing). Clearly, there are 
$u(n-1,k-1)$ such trees. On the other hand, if vertex 1 is not a leaf, 
amalgamate it with the root, erasing its parent edge and re-planting 
its child edges at the root. The pruned tree 
still has increasing leaves and its root has $j\ge k$ children. Also, 
the original tree can be recovered from the pruned tree provided
you know the position among the root edges of the leftmost re-planted edge (since the re-planted edges 
form a contiguous bunch of edges). There are $k$ choices for this 
position, independent of $j$, and so $k \sum_{j=k}^{n-1}u(n-1,j)$ choices altogether. Thus  
$\v\, \I_{n,k}\,\v$ does indeed satisfy recurrence (\ref{recur}).

\section{Stirling and Gessel permutations}\label{stirperm}

\vspace*{-3mm}

\hspace*{6mm} A size-$n$ Stirling permutation is a permutation of 
the multiset $\{1,1,2,2,\ldots,n,n\}$ in 
which, for each $i \in [n]$, all entries between the two occurrences of 
$i$ exceed $i$ (the Stirling property for $i$)\cite{stirpoly,dblfac}. Stirling permutations are 
counted by the odd double factorial, and there is a simple bijection from \IOTs 
to Stirling permutations due to Svante Janson \cite{janson08}:
given an \IOT, delete the root label and transfer the remaining labels from vertices to 
parent edges. Walk clockwise around the tree thereby traversing each 
edge twice and record labels in the order encountered; thus each 
label is recorded twice. \qed

As noted by Ira Gessel, permutations of 
the multiset $\{1,1,2,2,\ldots,n,n\}$ in 
which the second occurrences of $1,2,3,\ldots,n$ occur in that order, equivalently the second occurrence of each $i \in [n]$ is a right-to-left minimum 
(the Gessel property), are also counted by $(2n-1)!!$. Indeed, a 
Gessel permutation of size $n$ is obtained from one of size $n-1$ by 
inserting the first $n$ in any one of the $2n-1$ possible slots and 
the second second $n$ at the end. Similarly, a Stirling permutation of size $n$ 
is obtained from one of size $n-1$ by 
inserting the first $n$  in any one of the $2n-1$ possible slots and 
the second $n$ immediately after it. These observations lead immediately 
to a recursively-defined bijection $\phi $ from Gessel to Stirling 
permutations: 

\textbf{Recursive} $\mbf{\phi}$. \ 
Given a Gessel permutation $p$ of size $n$, let $j$ 
denote the position of the first $n$ in $p$, form $p'$ by deleting 
the two $n$'s from $p$; then $\phi(p)$ is obtained by inserting two 
$n$'s into the recursively defined $\phi(p')$ so that they occupy positions $j$ and $j+1$.
The base case is $\phi\big((1,1)\big)=(1,1)$. \qed

The bijection $\phi$ is worth a closer look because it involves, in a 
simpler setting, ideas used in our main bijection. It has 
two non-recursive descriptions.

\textbf{Direct algorithmic}  $\mbf{\phi}$. \  Suppose given a size-$n$ Gessel permutation. Find the largest $i$ 
whose two occurrences bracket an entry $<i$ (if there is no such 
such $i$, stop: the permutation is already Stirling).  Place a divider 
just after each occurrence of $i$ and circle (or ``ensquare'') each 
entry $\le i$ between the two dividers. Then cyclically shift right 
the contents of the squares. This step ensures the Stirling 
property holds for $i$, but at the expense of the Gessel property\,: 
the second occurrence of $i$  
now initiates an inversion in the permutation. Do likewise with the 
next largest $i$ whose two occurrences bracket an entry $<i$, and so 
on, until a Stirling permutation is obtained. 
For example, $p=2\ 3\ 1\ 5\ 4\ 1\ 2\ 3\ 4\ 6\
5\ 6$ is processed as follows:

\vspace*{-3mm}

{\scriptsize
\[
\begin{array}{cccccc}
  2\ 3\ 1\ 5\ 4\ 1\ 2\ 3\ 4\ 6\, 
\left|\!\vphantom{a_{a}^{A}}\ 
\framebox{5}\ \framebox{6}\, \right|  &  \rightarrow & 
2\ 3\ 1\ 5\ 4\ 1\ 2\ 3\ 4\ 6\ 
\gray1{\framebox{\black{6}}\ \framebox{\black{5}}} & = &
 \ 2\ 3\ 1\ 5\,
\left|\!\vphantom{a_{a}^{A}}\ 
\framebox{4}\ \framebox{1}\ \framebox{2}\ \framebox{3}\ \framebox{4}\ 6\ 6\ 
\framebox{5}\, \right|  & \rightarrow \\[3mm]
\ 2\ 3\ 1\ 5\ 
\gray1{\framebox{\black{5}}\ \framebox{\black{4}}\ \framebox{\black{1}}\ \framebox{\black{2}}\ \framebox{\black{3}}}\ 6\ 6\ 
\gray1{\framebox{\black{4}}}  &   = &   2\ 3\ 1\ 5\ 5\ 4\,
\left|\!\vphantom{a_{a}^{A}}\ 
\framebox{1}\ \framebox{2}\ \framebox{3}\ 6\ 6\ \framebox{4}\, 
\right| & \rightarrow & 2\ 3\ 1\ 5\ 5\ 4\ 
\gray1{\framebox{\black{4}}\ \framebox{\black{1}}\ \framebox{\black{2}}}\ 6\ 6\ \gray1{\framebox{\black{3}}} & =  \\[3mm]
2\ 3\  
\left|\!\vphantom{a_{a}^{A}}\ \framebox{1}\ 5\ 5\ 4\ 4\ \framebox{1}\ \framebox{2}\ 6\ 6\ 
\framebox{3}\, \right|  & \rightarrow &    
2\ 3\ \gray1{\framebox{\black{3}}}\ 5\ 5\ 4\ 4\ \gray1{\framebox{\black{1}}\ \framebox{\black{1}}}\ 6\ 6\ 
\gray1{\framebox{\black{2}}}  &  = & 
2\,\left|\!\vphantom{a_{a}^{A}}\ 
3\ 3\ 5\ 5\ 4\ 4\ \framebox{1}\ \framebox{1}\ 6\ 6\ 
\framebox{2}\, \right|  & \rightarrow  \\[3mm]
2\ 3\ 3\ 5\ 5\ 4\ 4\ \gray1{\framebox{\black{2}}\ \framebox{\black{1}}}\ 6\ 6\ 
\gray1{\framebox{\black{1}}}  & = &  2\ 3\ 3\ 5\ 5\ 4\ 4\ 2\ 1\ 6\ 6\ 
1   &   &   &     
\end{array}
\]
}

\vspace*{-3mm}

To reverse the map, start with the smallest $i$ whose second 
occurrence initiates an inversion. Take this second 
occurrence and all its inversion terminators, and rotate them left, 
and so on. \qed

\textbf{Constructive}  $\mbf{\phi}$ \textbf{via trapezoidal words}. \  Following Riordan, a \emph{trapezoidal word} 
$(w_{i})_{i=1}^{n}$ is a member of the Cartesian product
$[1]\times [3]\times \ldots \times [2n-1]$. Now $\phi$ is the 
composition of the following bijections.

\noindent \textbf{Gessel permutation $\mbf{\rightarrow}$ trapezoidal word}: Given a  Gessel permutation $\si$, for 
$1\le i \le n$ set $w_{i} =$ 
number of entries weakly preceding the first $i$ in $\si$ and $\le i$.

\noindent \textbf{Trapezoidal word $\mbf{\rightarrow}$ Gessel permutation}:  Given a  trapezoidal word 
$(w_{i})_{i=1}^{n}$, start with a row of $2n$ empty squares to be 
filled with the entries of the permutation. Place the first $n$ in 
the $w_{n}$-th square and the second $n$ in the last
square. Then for $i=n-1,n-2,\ldots,1$ in turn, place the first $i$ in 
the $w_{i}$-th unoccupied square and the second $i$ in the last
unoccupied square.

\noindent \textbf{Stirling permutation $\mbf{\rightarrow}$ trapezoidal word}: Given a  Stirling permutation $\si$, for 
$1\le i \le n$ set $w_{i} =$ 
number of entries weakly preceding the first $i$ in $\si$ and $\le i$.

\noindent \textbf{Trapezoidal word $\mbf{\rightarrow}$ Stirling permutation}:  Given a  trapezoidal word 
$(w_{i})_{i=1}^{n}$, start with a row of $2n$ empty squares to be 
filled with the entries of the permutation. Place the first $n$ in 
the $w_{n}$-th square and the second $n$ in the next
square. Then for $i=n-1,n-2,\ldots,1$ in turn, place the first $i$ in 
the $w_{i}$-th unoccupied square and the second $i$ in the next
unoccupied square.

For example, the Gessel permutation $p$ above gives the trapezoidal word $1\ 1\ 2\ 4\ 4\ 10$. 
Then $\phi(p)$ is constructed as follows (subscripts indicate 
placement positions among unoccupied squares).
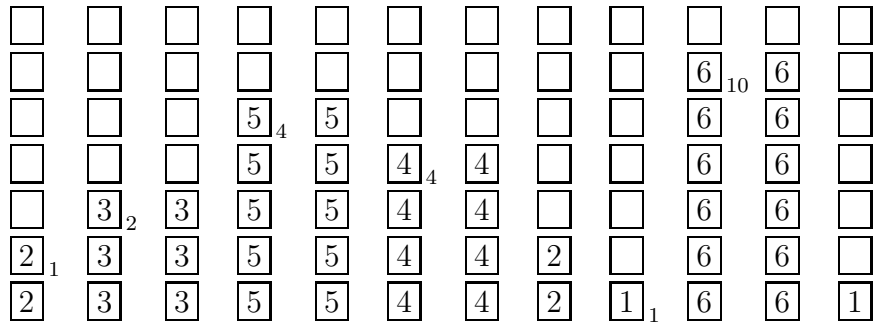
\[
\begin{array}{llllllllllll}
    \framebox{\phantom{2}} & \framebox{\phantom{2}} & 
    \framebox{\phantom{2}}_{\phantom{2}} & \framebox{\phantom{2}} & 
    \framebox{\phantom{2}}_{\phantom{2}} & \framebox{\phantom{2}} & 
    \framebox{\phantom{2}}_{\phantom{2}} & \framebox{\phantom{2}}_{\phantom{2}} & 
   \framebox{\phantom{2}}_{\phantom{2}} & \framebox{\phantom{2}} & 
   \framebox{\phantom{2}}_{\phantom{2}} & \framebox{\phantom{2}}  \\
   
    \framebox{\phantom{2}} & \framebox{\phantom{2}} & \framebox{\phantom{2}} & \framebox{\phantom{2}} & 
    \framebox{\phantom{2}} & \framebox{\phantom{2}} & \framebox{\phantom{2}} & \framebox{\phantom{2}} & 
    \framebox{\phantom{2}} & \framebox{6}_{\:10}\!\! & \framebox{6} & \framebox{\phantom{2}}  \\
    \framebox{\phantom{2}} & \framebox{\phantom{2}} & 
    \framebox{\phantom{2}} & \framebox{5}_{\:4} & 
    \framebox{5} & \framebox{\phantom{2}} & \framebox{\phantom{2}} & \framebox{\phantom{2}} & 
    \framebox{\phantom{2}} & \framebox{6} & \framebox{6} & \framebox{\phantom{2}}  \\
    \framebox{\phantom{2}} & \framebox{\phantom{2}} & \framebox{\phantom{2}} & \framebox{5} & 
    \framebox{5} & \framebox{4}_{\:4} & \framebox{4} & \framebox{\phantom{2}} & 
    \framebox{\phantom{2}} & \framebox{6} & \framebox{6} & \framebox{\phantom{2}}  \\
    \framebox{\phantom{2}} & \framebox{3}_{\:2} & \framebox{3} & \framebox{5} & 
    \framebox{5} & \framebox{4} & \framebox{4} & \framebox{\phantom{2}} & 
    \framebox{\phantom{2}} & \framebox{6} & \framebox{6} & \framebox{\phantom{2}}  \\
    \framebox{2}_{\:1} & \framebox{3} & \framebox{3} & \framebox{5} & 
    \framebox{5} & \framebox{4} & \framebox{4} & \framebox{2} & 
    \framebox{\phantom{2}} & \framebox{6} & \framebox{6} & \framebox{\phantom{2}} \\
    \framebox{2} & \framebox{3} & \framebox{3} & \framebox{5} & 
    \framebox{5} & \framebox{4} & \framebox{4} & \framebox{2} & 
    \framebox{1}_{\:1} & \framebox{6} & \framebox{6} & \framebox{1}
\end{array}
\]

\vspace*{-3mm}

\section{2-Configurations}\label{2config} 
Fix a positive integer $n$ and take a multiset $M$ consisting of the 
union of $[n]$ and an arbitrary subset of $n$. An element of $[n]$ 
that occurs twice in $M$ is called a \emph{repeater}. Consider a partition
of $M$ into a list of decreasing lists (blocks) such that 
\begin{enumerate}
    \vspace*{-3mm}
    \item  The entries in each block are distinct.

    \item  Each block contains at most one first occurrence of a repeater. 
    
    \item If a block contains a first occurrence of a repeater, then all later entries in the block are second occurrences of other repeaters.    
    
\end{enumerate}

Call such a partition a \emph{2-configuration} of size $n$. Draw an arc 
connecting the two occurrences of each repeater. We can then view the 
2-configuration (i) as an arc diagram, and speak of noncrossing arcs, (ii) as a graph  with the blocks as vertices, and speak of its (connected) components. By directing each arc from left to right we may also speak of incoming and outgoing arcs. In these terms, conditions 2 and 3 say that each block has at most one outgoing arc and all entries in a block after an outgoing arc have incoming arcs.

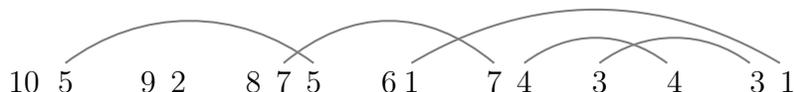
\begin{center}
\begin{pspicture}(-6,0)(4,0.8)
\psset{unit=1.0cm}

\rput(-5.85,0){\textrm{{\normalsize $10$}}}
\rput(-5.3,0){\textrm{{\normalsize $5$}}}
\rput(-4.2,0){\textrm{{\normalsize $9$}}}
\rput(-3.8,0){\textrm{{\normalsize $2$}}}
\rput(-2.8,0){\textrm{{\normalsize $8$}}}
\rput(-2.4,0){\textrm{{\normalsize $7$}}}
\rput(-2.0,0){\textrm{{\normalsize $5$}}}

\rput(-1.0,0){\textrm{{\normalsize $6$}}}
\rput(-.7,0){\textrm{{\normalsize $1$}}}

\rput(.4,0){\textrm{{\normalsize $7$}}}
\rput(.8,0){\textrm{{\normalsize $4$}}}


\rput(1.8,0){\textrm{{\normalsize $3$}}}

\rput(2.8,0){\textrm{{\normalsize $4$}}}

\rput(3.9,0){\textrm{{\normalsize $3$}}}
\rput(4.3,0){\textrm{{\normalsize $1$}}}

\psbezier[linecolor=gray,linewidth=.8pt](-5.3,.25)(-4.3,1.0)(-3,1.0)(-2,.25)

\psbezier[linecolor=gray,linewidth=.8pt](-2.4,.25)(-1.5,1.0)(-.5,1.0)(.4,.25)

\psbezier[linecolor=gray,linewidth=.8pt](.8,.25)(1.4,.7)(2.1,.7)(2.7,.25)

\psbezier[linecolor=gray,linewidth=.8pt](1.8,.25)(2.4,.7)(3.2,.7)(3.8,.25)
\psbezier[linecolor=gray,linewidth=.8pt](-.7,.25)(.9,1.2)(2.6,1.2)(4.2,.25)

\end{pspicture}
\end{center}
The example illustrated for $n=10$ has repeaters $1,\,3,\,4,\,5,\,7$; it has 8 
blocks and 5 arcs. It has 3 components, namely, $C_{1}=(6\ 1\, ,\ 3\, 
,\ 3\ 1),\ C_{2}=(9\ 2),\ C_{3}=(10\ 5\, ,8\ 7\ 5\, ,\ 7\ 
4\, ,\ 4)$. As here, the blocks in a component will always be listed in left-to-right order and the components themselves will always be ordered by increasing minimum entry. The first of these components has noncrossing arcs, the second is 
a singleton with one block and no arcs, and the third has some crossing arcs.  
For each nonsingleton component, its \emph{critical} block is the 
second block if the minimum entry in the component occurs in the first 
block, otherwise it is the first block. (Note this definition 
would not make sense for a singleton component.) Thus, 
the critical block for $C_{1}$ is the second one, (3), and for $C_{3}$ is the first one, (10\ 5).

A \emph{bad} repeater is a repeater whose two occurrences enclose a smaller 
entry in the 2-configuration. Any such smaller entry is a \emph{delinquent} 
for the bad repeater. In other words, a delinquent is an entry $i$ 
below an arc whose (equal) endpoints are $>i$. In the example, 
the bad repeaters are 5,\,7,\,4 with delinquent entries 
$\{2\},\ \{5,\,6,\,1\},\ \{3\}$ respectively.

Now we list some properties a 2-configuration may (or may not) have.
\begin{itemize}
    \vspace*{-3mm}
    \item  The \emph{good-component} property: \   For each nonsingleton 
    component $C_{i}$, all blocks in the components 
    $C_{1},C_{2},\ldots,C_{i-1}$ (including, of course, singleton components) precede the critical block of $C_{i}$ in 
    the 2-configuration. In the example, $C_{3}$ badly fails this condition.

    \item  The \emph{single-incoming-arc} property: \  
    No block has more than one incoming arc. In the example, 
    the last block fails this condition because it has 2 incoming arcs.

    \item  The \emph{no-crossing-in-component} property: \  The arcs within each 
    component, considered as connecting entries rather than blocks, are 
    noncrossing. In the example, $C_{3}$ fails this condition because 
    the arcs incident with the block $8\ 7\ 5$ do cross, while 
    $C_{1}$ and $C_{2}$ satisfy it.

    \item  The \emph{Gessel} property: \  The 
    second occurrence of each bad repeater is the last entry in its 
    block, and repeaters that end their second block are 
    right-to-left minima in the flattened 
    configuration (to flatten means to concatenate all the blocks). 
    In the example, the bad repeater 7 fails the first part of this 
    condition, and the repeaters 4,\,5 fail the second part.
    
     \item  The \emph{Stirling} property: \  No bad repeaters.
     
     \item  The \emph{restricted-first-entry} property: \  The first entry of a 
     block is never a repeater and these first entries increase left 
     to right.
\end{itemize}
\vspace*{-2mm}


\section{First steps in the main bijection}

\vspace*{-4mm}

\hspace*{6mm} The identification of (1-23-4)-avoiding permutations and 
\IOTs with increasing leaves can be described as a composition of four 
bijections:
\[
\I_{n} \longleftrightarrow \J_{n} \longleftrightarrow \s_{n} 
\longleftrightarrow \a_{n} \longleftrightarrow \p_{n},
\]
where $\I_{n}$ is the set of all size-$n$ \IOTs with 
increasing leaves, $\J_{n}$ is a certain subset of the size-$n$ Stirling 
permutations, $\s_{n}$ and $\a_{n}$ are 
subsets of the 2-configurations of size $n$, and $\p_{n}$ is the set of 
(1-23-4)-avoiding permutations of $[n]$.

\vspace*{-3mm}

\subsection{$\I_{n}\: \protect{\mbf{\rightarrow}}\: \J_{n}$}   

\vspace*{-3mm}

\hspace*{6mm}A \emph{plateau} in a Stirling permutation is a
pair of adjacent entries that are equal. Thus Janson's bijection of 
Section \ref{stirperm}
identifies leaves in the tree with plateaus in the permutation, and
sends $\I_{n}$, the size-$n$ \IOTs with 
increasing leaves, onto $\J_{n}$, defined as 
the set of size-$n$ Stirling permutations whose plateaus 
increase left to right. This is the first bijection.

\vspace*{-3mm}

\subsection[$\J_{n}\: \rightarrow\: \s_{n}$]{$\J_{n}\: \protect{\mbf{\rightarrow}}\: \s_{n}$}  

\vspace*{-3mm}

\hspace*{6mm}Given a permutation in $\J_{n}$, underline the longest 
decreasing run starting at the second entry of each plateau and 
extract these runs to form a set of decreasing lists (blocks) with 
increasing first entries:
\[
1\  3\  5\  \underline{5\  3\  1}\  2\  6\  7\  \underline{7\  6\  4}\  8\  9\ 
\underline{9\  8\  4\  2}\  10\  \underline{10}\ \rightarrow\ 5\: 3\: 
1 \quad\ 7\: 6\: 4 \quad\ 9\: 8\: 4\: 2 \quad\ 10
\]

We call the result $\rho$ a \emph{Stirling configuration}. 
The original permutation can be recovered from $\rho$. First, 
flatten $\rho$.
Then, for each $i$ that occurs only once, insert a second $i$ as far 
to the left of the existing $i$ as possible so that the two $i$'s 
don't enclose a smaller entry. (The order in which the missing $i$'s 
are inserted is immaterial.)

The first occurrence of a repeater in a Stirling 
configuration (4 is the only repeater in the example) is always the last entry in its block 
for otherwise a bad repeater would be present. 
The set of Stirling configurations of size $n$, denoted $\s_{n}$, 
can be characterized as the set of 2-configurations of size $n$ 
that have the restricted-first-entry and Stirling properties of 
the previous section.

Extraction of underlined runs is thus a bijection from 
$\J_{n}$, the increasing-plateau Stirling permutations, to 
$\s_{n}$. This is the second bijection. 

\vspace*{-3mm}

\subsection{$\p_{n}\: \mbf{\rightarrow}\: \a_{n}$} \label{PtoA}

\vspace*{-3mm}

\hspace*{6mm}Now we turn to the pattern-avoiding permutations. 
An \emph{ascent} in a permutation 
is a pair of adjacent entries $a,b$ with $a<b$\,; $a$ is the 
\emph{initiator}, $b$ the \emph{terminator}. Clearly, a permutation 
is (1-23-4)-avoiding if and only if, for each ascent, either its 
initiator is a left-to-right minimum or its terminator is a 
right-to-left maximum. But an asymmetrical viewpoint is more 
convenient.
Any permutation can be split into  
\emph{LRMin segments} starting at its left-to-right minima. The first two entries of a nonsingleton segment necessarily form an ascent, an \emph{LRMin ascent} (for otherwise the second entry would start a new segment);  other ascents are \emph{free}.
Then a permutation is 
\vspace*{-3mm}
\begin{itemize}
\item  (1-23)-avoiding iff it has no free ascents

\item  (1-23-4)-avoiding iff each free ascent terminates at a 
right-to-left maximum.
\end{itemize}


Suppose given a (1-23-4)-avoiding permutation. We will perform a 
series of six reversible steps, using
\[
23\ \,4\ \,21\ \,6\ \,25\ \,24\ \,14\ \,22\ \,18\ \,20\ \,16\ \,13\ \,11\ \, 19 \ \,7\ \,5\ \,  2\ \,8\ \,17\ \,12\ \,10\ \,  1\ \,15\ \,3\ \,9
\]
as a working example. The resulting configurations will form $\a_{n}$.

\noindent Step 1. Split the permutation into 
its LRMin segments. Working left to right, for each 
ascent initiator $i$, place an overline starting at $i$ and extending as far right 
as possible so that it does not cover an entry $<i$. Each new overline is placed below any existing overlines. Also, if a segment has just one entry and thus no ascents, overline that entry. Now every entry is covered by a unique overline (``cover'' means no intervening overlines). Clearly,  no overline will straddle two segments and 
each pair of overlines is either disjoint or nested (no overlaps). Also, an overline starting at a free ascent initiator will cover both entries of the ascent.
\[
\begin{array}{rcccclc}
  \gray1{\overline{\rule{0ex}{2ex}\black{23}}} \!\!\!&\textrm{ \Large \textbf{/}} &  \ \   \gray1{\overline{\black{4\ 21}\ \gray1{\overline{\black{6\ 25\ 24}\ \overline{\black{14\ 
   22}\   \overline{\rule{0ex}{2ex}\black{18\ 20}} \         \black{16}}\ 
   \black{13}\ \overline{\rule{0ex}{2ex}\black{11\ 19}}\ \black{7}}}\ \black{5}}} \!\!& 
\textrm{ \Large \textbf{/}} & \  \gray1{\overline{\rule{0ex}{2.5ex}\black{2}\  \gray1{\overline{\rule{0ex}{2ex}\black{8\ 
   17\ 12\ 10}}}}} \!\!\!\! & \ \    \textrm{ \Large \textbf{/}}  \ \
   &  \gray1{\overline{\rule{0ex}{2ex}\black{1\ 15}\ \gray1{\overline{\rule{0ex}{2ex}\black{3\ 9}}}}}  
\end{array}
\]

\noindent Step 2. 
For each LRMin segment,  extract the entries covered by an
overline to obtain a list of blocks. (With a slight abuse of  words, we will refer to this list as a segment of blocks.)
\[
23\!  \ \textrm{ \Large \textbf{/}}  \  \ 4\ 21\ 5  \quad 6\ 25\ 24\ 13\ 7   \quad 14\ 22\ 16 \quad 18\ 20 
\quad 11\ 19   \   \textrm{ \Large \textbf{/}}   \  \  2 \quad 8\  17\ 12\ 
10    \  \textrm{ \Large \textbf{/}}  \  \     1\ 15 
\quad\, 3\ 9
\]

To recapture the original LRMin segments,  ``coalesce''  the blocks in each segment, inserting each one (starting at the second block) into its current 
predecessor just far enough from 
the end so that its first entry does not terminate a free ascent. 


\noindent Step 3. Only the first block in a segment can be a singleton block. Furthermore, if the first block is a singleton and and the segment has more than one block, then the second entry of the second block, being a free ascent terminator, is a right-to-left maximum and, in particular, the largest entry in the segment. In this case, interchange the first and third (smallest and largest) entries in the segment. In the example step 3 alters the third segment.
\[
23\  \ \textrm{ \Large \textbf{/}}  \  \ 4\ 21\ 5  \quad 6\ 25\ 24\ 13\ 7   \quad 14\ 22\ 16 \quad 18\ 20 
\quad 11\ 19 \quad   \textrm{ \Large \textbf{/}}  \quad 17 \quad 8\  2\ 12\ 
10  \quad   \textrm{ \Large \textbf{/}}  \quad    1\ 15 
\quad\, 3\ 9
\]

\noindent Step 4. 
Sort each block into decreasing order.
\[
23\  \ \textrm{ \Large \textbf{/}}  \  \ 21\ 5\ 4  \quad 25\ 24\ 13\ 7 \ 6  \quad 22\ 16\ 14 \quad 20 \ 18
\quad 19\ 11 \quad   \textrm{ \Large \textbf{/}}  \quad 12\ 
10\ 8\ 2 \quad 17  \  \ ,   \textrm{ \Large \textbf{/}}  \quad    
15\ 1 \quad 9\ 3
\]
This step amounts to transferring the first entry of each block to 
the end of the block except, in case the first block in a segment is a 
singleton, transferring the first two entries of the second block (if there is one) to 
the end of the block. Thus the step can be reversed.

\noindent Step 5. Sort blocks in each segment in order of increasing 
first entries.
\[
23\  \ \textrm{ \Large \textbf{/}}  \  \ 19\ 11  \quad  20\ 18  \quad 21\ 5\ 4 \quad  22\ 16\ 14  \quad 25\ 
24\ 13\ 7\ 6 \quad   \textrm{ \Large \textbf{/}}  \quad 12\ 
10\ 8\ 2 \quad 17  \    \textrm{ \Large \textbf{/}}  \quad    
 9\ 3 \quad  15\ \,1 
\]
Originally, the block containing the smallest entry of a segment was 
the first block in the segment and the remaining blocks were in 
order of decreasing first entries, both statements holding unless the telltale singleton block is present in which case \emph{all} blocks were in order of decreasing first entries. So this step is reversible.

\noindent Step 6. Establish links between the blocks in a segment by inserting a second copy of the last entry of each 
block (save for the last one) into the next block, maintaining monotonicity.
\[
\begin{array}{l}
23\  \ \textrm{ \Large \textbf{/}}  \  \ 
    19\ 11  \quad  20\ 18\ 11  \quad 21\ 18\ 5\ 4 \quad  22\ 16\ 14\ 4  \quad 25\ 
24\ 14\ 13\ 7\ 6 \quad   \textrm{ \Large \textbf{/}}  \\[3mm] 
12\ 10\ 8\ 2 \quad 17\ 2  \quad   \textrm{ \Large \textbf{/}}  
 \quad   9\ 3 \quad  15\ 3\ 1 
\end{array}
\]
This step is obviously reversible.

\noindent Step 7. Erase dividers and sort the entire collection of blocks in order of increasing first entry.
\[
9\ 3 \quad  12\ 10\ 8\ 2     \quad  15\ 3\ 1  \quad 17\ 
 2 \quad  19\ 11  \quad  20\ 18\ 11  \quad 21\ 18\ 5\ 4 \quad  22\ 16\ 14\ 4 \quad 23   \quad 25\ 
24\ 14\ 13\ 7\ 6 
 \]
To reverse Step 7, follow the trail of the repeaters to determine the blocks in each segment. The order of blocks within a segment is determined by their first entries, which 
must be increasing, and the 
order of the segments themselves is determined by their smallest entries, which 
must be decreasing.

We call the result of Step 7 an \emph{avoider configuration}.  Avoider 
configurations form the set $\a_{n}$ and they are 
characterized as 2-configurations that have the restricted-first-entry and 
single-arc properties (both by construction) and also the good-component  
property (reflecting the fact that each free ascent in the permutation 
terminates at a right-to-left maximum).

This is the fourth of the four bijections. In the next section we 
present a bijection from the 2-configurations of size $n$ with the  
single-arc and good-component properties to 
the 2-configurations of size $n$ with the  
Stirling property. The restriction of this bijection to 2-configurations with the restricted-first-entry property then completes the proof, by providing the third of the desired bijections, from $\a_{n}$ to $\s_{n}$.

\vspace*{-3mm}

\section{Bijections between sets of 2-configurations}

\vspace*{-3mm}

\hspace*{6mm}The following four subsets of the 2-configurations of size $n$ turn out to be equinumerous.
\begin{enumerate}
    \vspace*{-2mm}
    \item  $\w_{n}:$ \ The 2-configurations that have both the single-incoming-arc and the
    good-component property. 

    \item  $\x_{n}:$ \ The 2-configurations that have both the 
    good-component and the no-crossing-in-component property.

    \item  $\y_{n}:$ \ The 2-configurations that have both the 
    no-crossing-in-component and the Gessel property.

    \item  $\z_{n}:$ \ The 2-configurations that have the 
    Stirling property.
\end{enumerate}

\vspace*{-2mm}

\noindent We will give bijections $\w_{n} \rightarrow  \x_{n} \rightarrow  \y_{n}
\rightarrow  \z_{n}$. A 2-configuration in $\z_{n}$ has no bad repeaters, so we must turn bad repeaters in $\sigma \in \w_{n}$ into good ones. Say a bad repeater is of Type 1 if its first occurrence is not the last entry in its block, of Type 2 if its second occurrence is 
not the last entry in its block, and of Type 3 if it is neither of 
Type 1 or 2. (It may be both Type 1 and 2). The three bijections eliminate the bad repeaters of Type 1,\,2, and 3 in turn. All three bijections  preserve the statistic ``number of components''. The first also preserves the support multiset of each component while the other two do not.

\noindent \textbf{Iterative descriptions:}
\vspace*{-5mm}
 
\begin{enumerate}
\item    $\w_{n} \rightarrow  \x_{n}$ \quad As long as there is a bad repeater of Type 1 in the current 2-configuration, find the largest one and transfer the last entry of its first block to its second block.

\item   $\x_{n} \rightarrow  \y_{n}$ \quad As long as there is a bad repeater of Type 2 in the current 2-configuration, find the largest one and transfer its second occurrence to the block of its leftmost delinquent. 

\item   $\y_{n} \rightarrow  \z_{n}$\quad As long as there is a bad repeater of Type 3in the current 2-configuration, find the largest one. Let 
$A_{1},A_{2},\ldots,A_{k}$ ($k\ge 1$) denote the blocks containing its 
delinquents and let $A_{k+1}$ denote 
the second block of the bad repeater. Transfer the bad repeater from $A_{k+1}$ 
to $A_{1}$ and transfer 
all delinquents in $A_{i}$ to $A_{i+1},\  1\le i \le k$. The illustration shows bad repeaters in \blue{blue}, delinquents in \red{red}. 
\[
\begin{array}{c}
3 \hspace*{10mm} \blue{7} \hspace*{10mm} 8\ \red{5} \hspace*{10mm} 9\ \red{4\ 2} \hspace*{10mm} \red{1} \hspace*{10mm} \red{6\ 3} \hspace*{10mm} \blue{7} \quad \rightarrow \\[1mm]
\blue{3}\hspace*{10mm}  7\hspace*{10mm}  8\ 7 \hspace*{10mm}  9\ 5 \hspace*{10mm}  4\ \red{2} \hspace*{10mm}  \red{1} \hspace*{10mm}  6\ \blue{3} \quad \rightarrow \\[1mm]  
3 \hspace*{10mm}  7 \hspace*{10mm}  8\ 7 \hspace*{10mm}  9\ 5 \hspace*{10mm}  4\ 3 \hspace*{10mm}  2 \hspace*{10mm}  6\ 1 \phantom{\ \quad \rightarrow} \\[3mm]
\textrm{\small The bijection $\y_{n} \rightarrow  \z_{n}$\hspace*{10mm}} 
\end{array}
\]
\end{enumerate}

\noindent \textbf{More direct descriptions:}
\vspace*{-5mm}

\begin{enumerate}
\item  $\w_{n} \rightarrow  \x_{n}$ \quad Given $\sigma \in \w_{n}$, transform each component as illustrated.

\begin{center}
\begin{pspicture}(-11,-2.6)(1,.6)
\psset{unit=1.0cm}

\rput(-10.5,0){\textrm{{\normalsize $13\ 5$}}}
\rput(-9,0){\textrm{{\normalsize $8\ 5$}}}
\rput(-7.6,0){\textrm{{\normalsize $8\ 3$}}}
\rput(-5.8,0){\textrm{{\normalsize $11\ 4\ \ 3$}}}
\rput(-4,0){\textrm{{\normalsize $9\ \,4$}}}
\rput(-2.5,0){\textrm{{\normalsize $9\ 6$}}}
\rput(-0.9,0){\textrm{{\normalsize $6\ \ 2$}}}
\rput(0.4,0){\textrm{{\normalsize $2$}}}

\psbezier[linecolor=gray,linewidth=.8pt](-10.25,.25)(-10,.6)(-9.1,.6)(-8.85,.25)
\psbezier[linecolor=gray,linewidth=.8pt](-9.15,.25)(-8.9,.6)(-8.0,.6)(-7.75,.25)
\psbezier[linecolor=gray,linewidth=.8pt](-7.4,.25)(-6.9,.7)(-5.8,.7)(-5.3,.25)
\psbezier[linecolor=gray,linewidth=.8pt](-5.7,.25)(-5.3,.6)(-4.2,.6)(-3.8,.25)
\psbezier[linecolor=gray,linewidth=.8pt](-4.15,.25)(-3.75,.6)(-3.1,.6)(-2.7,.25)
\psbezier[linecolor=gray,linewidth=.8pt](-2.3,.25)(-2.0,.6)(-1.5,.6)(-1.2,.25)
\psbezier[linecolor=gray,linewidth=.8pt](-0.7,.25)(-0.4,.6)(0.1,.6)(0.4,.25)

\psbezier[linecolor=blue,linewidth=.8pt,arrows=->,arrowsize=4pt 
3](-8.8,-.25)(-8.5,-.6)(-7.9,-.6)(-7.6,-.25)
\psbezier[linecolor=blue,linewidth=.8pt,arrows=->,arrowsize=4pt 
3](-5.3,-.25)(-4.9,-1)(-1.2,-1)(-0.8,-.25)
\psbezier[linecolor=blue,linewidth=.8pt,arrows=->,arrowsize=4pt 
3](-3.7,-.25)(-3.4,-0.7)(-1.3,-0.7)(-1,-.25)

\rput(-5.6,-1.6){\textrm{{\normalsize $\downarrow$}}}

\end{pspicture}
\end{center}

\begin{center}
\begin{pspicture}(-6.5,-1.0)(5.5,0.5)
\psset{unit=1.0cm}

\rput(-6.5,0){\textrm{{\normalsize $13\ 5$}}}
\rput(-5,0){\textrm{{\normalsize $8$}}}
\rput(-3.2,0){\textrm{{\normalsize $8\ 5\ 3$}}}
\rput(-1.3,0){\textrm{{\normalsize $11\ 4$}}}
\rput(0.1,0){\textrm{{\normalsize $9$}}}
\rput(1.5,0){\textrm{{\normalsize $9\ 6$}}}
\rput(3.5,0){\textrm{{\normalsize $6\ 4\ 3\ 2$}}}
\rput(5.2,0){\textrm{{\normalsize $2$}}}

\psbezier[linecolor=gray,linewidth=.8pt](-6.25,.25)(-5.85,.8)(-3.6,.8)(-3.2,.25)
\psbezier[linecolor=gray,linewidth=.8pt](-5,.25)(-4.6,.5)(-3.95,.5)(-3.55,.25)
\psbezier[linecolor=gray,linewidth=.8pt](-2.85,.25)(-2.4,1.2)(3.2,1.2)(3.65,.25)
\psbezier[linecolor=gray,linewidth=.8pt](0.1,.25)(.3,.5)(1.1,.5)(1.3,.25)
\psbezier[linecolor=gray,linewidth=.8pt](-1.05,.25)(-0.55,.9)(2.85,.9)(3.35,.25)
\psbezier[linecolor=gray,linewidth=.8pt](1.7,.25)(2.1,.5)(2.6,.5)(3.0,.25)
\psbezier[linecolor=gray,linewidth=.8pt](4,.25)(4.3,.5)(4.9,.5)(5.2,.25)

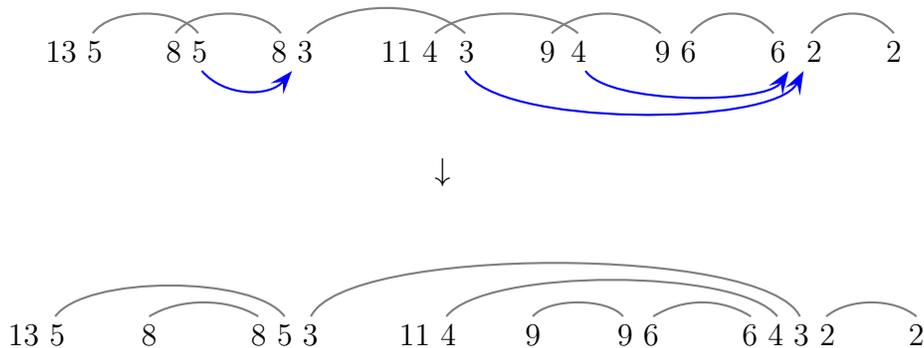
\rput(-1,-1.0){\textrm{{\small A component and its image under the bijection $\w_{n} \rightarrow  \x_{n}$}}}

\end{pspicture}
\end{center}
For each bad repeater $i$ of Type 1, transfer its second occurrence rightward to the first later block in its component whose outgoing arc starts at an entry $<i$, or to the last block in the component if there is no such block. The single-incoming-arc property has been traded for the no-crossing-in-component property, and the resulting configuration has no bad repeaters of Type 1.

\item  $\x_{n} \rightarrow  \y_{n}$ \quad Transfer the second occurrence of each bad repeater $i$ of Type 2 to the block of its leftmost delinquent, as illustrated by the blue arrows. For the reader's convenience, the components are indicated: $C_{i}$ is placed directly above the $i$-th component's unique block if the $i$-th component is a singleton, and above its critical block otherwise. (The color coding is merely for the reader's convenience.)

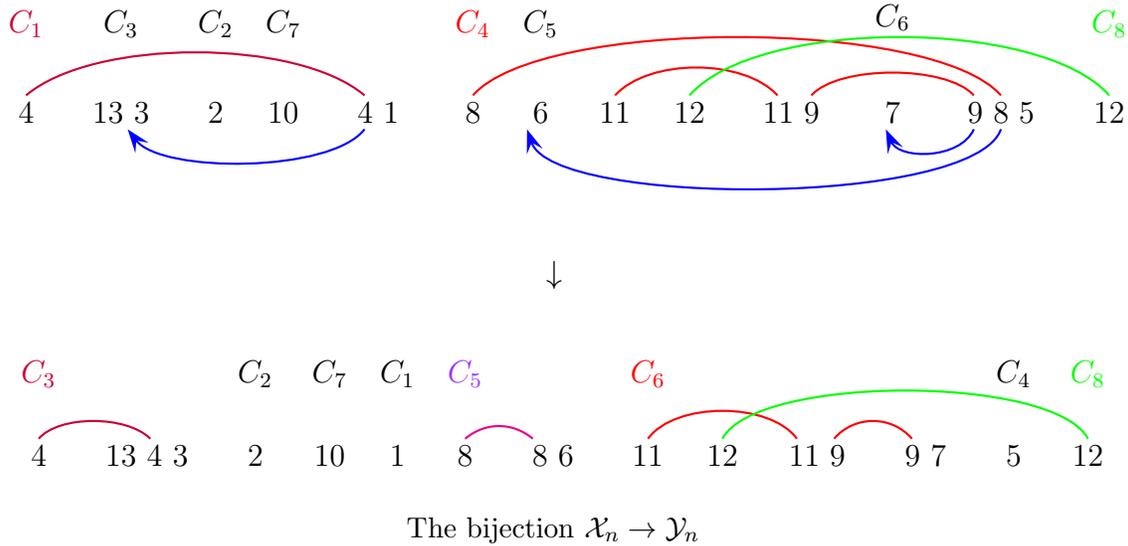
\begin{center}
\begin{pspicture}(-11.5,-2.7)(4.5,2)
\psset{unit=0.9cm}

\rput(-12.8,0){\textrm{{\normalsize $4$}}}
\rput(-11.4,0){\textrm{{\normalsize $13\ 3$}}}
\rput(-10,0){\textrm{{\normalsize $2$}}}
\rput(-9,0){\textrm{{\normalsize $10$}}}
\rput(-7.6,0){\textrm{{\normalsize $4\ 1$}}}
\rput(-6.2,0){\textrm{{\normalsize $8$}}}
\rput(-5.2,0){\textrm{{\normalsize $6$}}}
\rput(-4.1,0){\textrm{{\normalsize $11$}}}
\rput(-3.0,0){\textrm{{\normalsize $12$}}}
\rput(-1.5,0){\textrm{{\normalsize $11\ 9$}}}
\rput(0,0){\textrm{{\normalsize $7$}}}
\rput(1.6,0){\textrm{{\normalsize $9\ 8\ 5$}}}
\rput(3.2,0){\textrm{{\normalsize $12$}}}

\rput(-12.8,1.3){\textrm{{\normalsize \browne{$C_{1}$}}}}
\rput(-11.4,1.3){\textrm{{\normalsize $C_{3}$}}}
\rput(-10,1.3){\textrm{{\normalsize $C_{2}$}}}
\rput(-9,1.3){\textrm{{\normalsize $C_{7}$}}}

\rput(-6.2,1.3){\textrm{{\normalsize \red{$C_{4}$}}}}
\rput(-5.2,1.3){\textrm{{\normalsize $C_{5}$}}}

\rput(0,1.4){\textrm{{\normalsize $C_{6}$}}}

\rput(3.2,1.3){\textrm{{\normalsize \green{$C_{8}$}}}}

\psbezier[linecolor=brown,linewidth=.8pt](-12.8,.25)(-11.8,1.1)(-8.8,1.1)(-7.8,.25)
\psbezier[linecolor=red,linewidth=.8pt](-4.1,.25)(-3.6,.8)(-2.2,.8)(-1.7,.25)
\psbezier[linecolor=red,linewidth=.8pt](-6.2,.25)(-5.2,1.4)(0.6,1.4)(1.6,.25)
\psbezier[linecolor=red,linewidth=.8pt](-1.2,.25)(-.8,.7)(.8,.7)(1.2,.25)
\psbezier[linecolor=green,linewidth=.8pt](-3,.25)(-2,1.4)(2.2,1.4)(3.2,.25)

\psbezier[linecolor=blue,linewidth=.8pt,arrows=->,arrowsize=4pt 
3](-7.8,-.25)(-8.3,-.9)(-10.8,-.9)(-11.3,-.25)

\psbezier[linecolor=blue,linewidth=.8pt,arrows=->,arrowsize=4pt 
3](1.2,-.25)(1,-.7)(0.1,-.7)(-0.1,-.25)

\psbezier[linecolor=blue,linewidth=.8pt,arrows=->,arrowsize=4pt 
3](1.6,-.25)(1.2,-1.4)(-5.0,-1.4)(-5.4,-.25)

\rput(-5,-2.4){\textrm{{\normalsize $\downarrow$}}}

\end{pspicture}
\end{center}
\begin{center}
\begin{pspicture}(-6.5,-1.3)(7.5,1.2)
\psset{unit=0.9cm}

\rput(-7.6,0){\textrm{{\normalsize $4$}}}
\rput(-6,0){\textrm{{\normalsize $13\ 4\ 3$}}}
\rput(-4.4,0){\textrm{{\normalsize $2$}}}
\rput(-3.3,0){\textrm{{\normalsize $10$}}}
\rput(-2.3,0){\textrm{{\normalsize $1$}}}
\rput(-1.3,0){\textrm{{\normalsize $8$}}}
\rput(0,0){\textrm{{\normalsize $8\ 6$}}}
\rput(1.4,0){\textrm{{\normalsize $11$}}}
\rput(2.5,0){\textrm{{\normalsize $12$}}}
\rput(3.9,0){\textrm{{\normalsize $11\ 9$}}}
\rput(5.5,0){\textrm{{\normalsize $9\ 7$}}}
\rput(6.8,0){\textrm{{\normalsize $5$}}}
\rput(7.9,0){\textrm{{\normalsize $12$}}}

\psbezier[linecolor=brown,linewidth=.8pt](-7.6,.25)(-7.3,.6)(-6.3,.6)(-5.95,.25)
\psbezier[linecolor=magenta,linewidth=.8pt](-1.3,.25)(-1,.5)(-0.6,.5)(-0.3,.25)

\psbezier[linecolor=red,linewidth=.8pt](1.4,.25)(1.8,.8)(3.2,.8)(3.6,.25)
\psbezier[linecolor=red,linewidth=.8pt](4.15,.25)(4.45,.6)(5.0,.6)(5.3,.25)
\psbezier[linecolor=green,linewidth=.8pt](2.5,.25)(3,1.2)(7.4,1.2)(7.9,.25)

\rput(-7.6,1.2){\textrm{{\normalsize \browne{$C_{3}$}}}}

\rput(-4.4,1.2){\textrm{{\normalsize $C_{2}$}}}
\rput(-3.3,1.2){\textrm{{\normalsize $C_{7}$}}}
\rput(-2.3,1.2){\textrm{{\normalsize $C_{1}$}}}
\rput(-1.3,1.2){\textrm{{\normalsize \maggy{$C_{5}$}}}}

\rput(1.4,1.2){\textrm{{\normalsize \red{$C_{6}$}}}}

\rput(6.8,1.2){\textrm{{\normalsize $C_{4}$}}}
\rput(7.9,1.2){\textrm{{\normalsize \green{$C_{8}$}}}}

\rput(0,-1.1){\textrm{{\small The bijection $\x_{n} \rightarrow  \y_{n}$}}}

\end{pspicture}
\end{center}
The good-component property has been traded for the Gessel property.

\item  $\y_{n} \rightarrow  \z_{n}$\quad  
For each $i\in [n]$, record the
number $w_{i}$ of blocks that weakly precede the first appearance of $i$ 
and contain an entry $\le i$. 
Define an \emph{ender} to be an element of 
$[n]$ that appears as the last entry in a block.
Start with the same number of blocks, initially empty, as in the given 2-configuration. 
For each ender $i$ from largest to smallest in turn, place a copy of 
$i$ in the $w_{i}$-th empty block and, if $i$ is a repeater, place a 
second copy in the next empty block. When done, each block contains 
one entry. Next, for each non-ender $i$, place $i$ in the $w_{i}$-th 
\emph{available} block where the available blocks for $i$ are those that 
don't contain an ender $>i$. (The order in which the non-enders are 
placed is immaterial.) The resulting configuration is in $\z_{n}$. For 
the example above, the enders are $1,2,3,5,7$ of 
which $3,7$ are repeaters, and the $w_{i}$ are as follows.
\begin{center}
\begin{tabular}{|c||c|c|c|c|c||c|c|c|c|}\hline
     &  \multicolumn{5}{c||}{\rule[-3mm]{0mm}{8mm}enders}& 
     \multicolumn{4}{c|}{non-enders}  \\ \hline
    $i$ & 1 & 2 & 3 & 5 & 7 & 4 & 6 & 8 & 9  \\
    $w_{i}$ & 1 & 1 & 1 & 2 & 2 & 2 & 5 & 3 & 4  \\ \hline
\end{tabular}
\end{center}
The enders are placed successively as
illustrated in the following table.
\[
\begin{array}{|r|r|r|r|r|r|r|}\hline
    \phantom{1234} &\phantom{123}7  & \phantom{123}7 & \phantom{1234} & 
    \phantom{1234} & \phantom{1234} & \phantom{1234}  \\ \hline
     & 7 & 7 & 5 &  &  &   \\ \hline
    3 & 7 & 7 & 5 & 3 &  &   \\ \hline
    3 & 7 & 7 & 5 & 3 & 2 &   \\ \hline
    3 & 7 & 7 & 5 & 3 & 2 & 1 \\ \hline
\end{array}
\]

\noindent As for the non-enders, the available blocks for $i=4$ are the first, 
fourth, etc., and since $w_{i}=2$, 4 goes in the fourth block. For 
$i=6$ only the second and third blocks are unavailable and  
$w_{i}=5$ puts 6 in the last block. All blocks are available for 8 
and 9 and the final result is the same as in the first description. \qed

\end{enumerate}

\vspace*{-2mm}

It is fairly easy to see that the steps are reversible for all three mappings and that they do map onto the claimed sets.

\section{Another manifestation of the counting sequence}

\vspace*{-4mm}

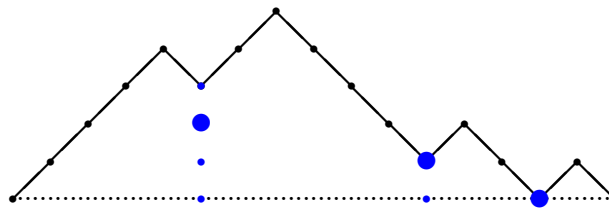
\hspace*{6mm}A \emph{valley-marked} Dyck path is a Dyck path in which, for each valley 
($DU$), one of the lattice points between the valley vertex and the 
$x$-axis inclusive (in blue in the figure below) is marked.
\Einheit=0.5cm
\[
\Pfad(-8,0),3333433444434434\endPfad
\SPfad(-8,0),1111111111111111\endSPfad
\DuennPunkt(-8,0)
\DuennPunkt(-7,1)
\DuennPunkt(-6,2)
\DuennPunkt(-5,3)
\DuennPunkt(-4,4)
\DuennPunkt(-3,3)
\DuennPunkt(-2,4)
\DuennPunkt(-1,5)
\DuennPunkt(0,4)
\DuennPunkt(1,3)
\DuennPunkt(2,2)
\DuennPunkt(3,1)
\DuennPunkt(4,2)
\DuennPunkt(5,1)
\DuennPunkt(6,0)
\DuennPunkt(7,1)
\DuennPunkt(8,0)
\blue{\DuennPunkt(3,0)
\DuennPunkt(-3,0)
\DuennPunkt(-3,1)
\DuennPunkt(-3,3)
\DickPunkt(-3,2)
\DickPunkt(3,1)
\DickPunkt(6,0)}
\Label\u{\textrm{A valley-marked Dyck path}}(0,-1)
\]

We have seen in Section \ref{iot} that recurrence (\ref{recur}) counts 
$\I_{n,k}$. It also counts $\d_{n,k}$, the valley-marked Dyck paths of 
semilength $n$ with first ascent of length $k$. To see this, 
observe that deleting the first peak produces a smaller-size valley-marked Dyck path. 
So every path in $\d_{n,k}$ 
is obtained by inserting a new first peak at height $k$ either into a path 
in $\d_{n-1,k-1}$ (no new valley is created) or into a path 
in $\d_{n-1,j}$ with $j\ge k$ (a new valley is created with $k$ 
choices for its marked lattice point), and the recurrence follows. From their common recurrence it is 
easy to construct a recursively-defined bijection between $\d_{n,k}$ 
and $\I_{n,k}$. But we can also give a direct bijection by exhibiting 
identical codings for $\d_{n,k}$ and 
$\I_{n,k}$. Given a marked path in $\d_{n,k}$, let $a_{i}=\#\,U$s between 
the $i$-th $D$ and the $(i+1)$-st $D$, $0\le i\le n-1$. 
Let $h_{i}$ denote the height above the 
line $y=-1$ of the mark for the valley corresponding to each
$i\in[1,n-1]$ with $a_{i}\ge 1$. The $a_{i}$'s code the path and the 
$h_{i}$'s code the marks. The example above yields
\[
\begin{array}{rcccccccc}
    i= & 0 & 1 & 2 & 3 & 4 & 5 & 6 & 7  \\
    a_{i}= & 4 & 2 & 0 & 0 & 0 & 1 & 0 & 1  \\
    h_{i}= &  & 3 &  &  &  & 2 &  & 1.
\end{array}
\]
For an ordered tree, to prune a non-root vertex means to delete the 
vertex and its parent edge and attach all its children to its parent, 
maintaining order. Now given a tree in $\I_{n,k}$, let $a_{i}$ denote 
the outdegree (= \#\,children) of vertex $i,\ 0 \le i \le n-1$ (vertex $n$ 
is necessarily a leaf), and for each non-root non-leaf vertex $i$, 
let $h_{i}$ denote the position of $i$ among its siblings in the tree 
obtained by pruning vertices $1,2,\ldots,i-1$. For example, the tree 
below yields the same code as the path above.

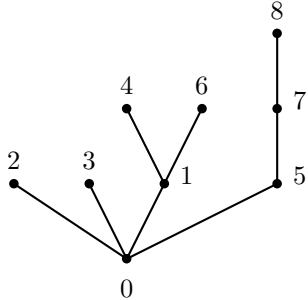
\begin{center} 

\begin{pspicture}(-7,-0.5)(9,3.5)
\psset{unit=1cm}

\psline(-1.5,1)(0,0)(-.5,1)
\psline(0,2)(.5,1)(1,2)
\psline(.5,1)(0,0)(2,1)(2,2)(2,3)

\psdots(-1.5,1)(-.5,1)(.5,1)(2,1)(0,0)(0,2)(1,2)(2,2)(2,3)

\rput(0,-0.4){\textrm{{\footnotesize $0$}}}
\rput(-1.5,1.3){\textrm{{\footnotesize $2$}}}
\rput(-.5,1.3){\textrm{{\footnotesize $3$}}}
\rput(.8,1.1){\textrm{{\footnotesize $1$}}}
\rput(2.3,1.1){\textrm{{\footnotesize $5$}}}
\rput(2.3,2.1){\textrm{{\footnotesize $7$}}}
\rput(2,3.3){\textrm{{\footnotesize $8$}}}
\rput(0,2.3){\textrm{{\footnotesize $4$}}}
\rput(1,2.3){\textrm{{\footnotesize $6$}}}

\end{pspicture}
\end{center}

We leave the reader to verify that both paths and trees can be 
uniquely retrieved from their codes, and that the codings are identical: 
pairs of sequences of nonnegative integers 
$\big( (a_{i})_{i=0}^{n-1},(h_{i})_{i\ge 1,\,a_{i}\ge 1}\big)$ such 
that $a_{0}\ge 1,\ a_{0}+a_{1}\ge 2,\ldots, a_{0}+a_{1}+\ldots +
a_{n-1} \ge n$ with equality in the last inequality, and $1\le 
h_{i} \le \sum_{j=0}^{i-1}a_{j}-(i-1)$ for each $i\ge 1$ such that 
$a_{i}\ge 1$.

\section{A generating function}

\vspace*{-4mm}

\hspace*{6mm} In  Section \ref{PtoA} we considered the LRMin segments of a permutation and distinguished between LRMin and free ascents. Call a (LRMin) segment \emph{short} if it has length 1, otherwise \emph{long}. There is an elegant generating function to count all of these notions for general permutations. Let $u(n,i,j,k)$ denote the number of permutations of $[n]$ with $i$ short segments, $j$ long segments (and hence $j$ LRMin ascents), and $k$ free ascents.
It is straightforward to obtain the recurrence 
\begin{multline*}
u(n,i,j,k)=u(n-1,i-1,j,k) + (j+k)u(n-1,i,j,k) + \\
(i+1)u(n-1,i+1,j-1,k) + (n-i-j-k)u(n-1,i,j,k-1),
\end{multline*}
valid for $n\ge 2,\:i\ge 1,\:j\ge 0,\:k\ge 0,\:i+j\ge 1,\:i+2j\le n,\:k\le n-i-j$, with initial condition $u(1,1,0,0)=1$.

Set $F(x,y,z,w):=1+\sum_{n\ge 1,i\ge 1,j\ge 0,k\ge 0}u(n,i,j,k)x^{n}/n!\:y^{i} z^{j} w^{k}$. 
The recurrence for $u$ translates into a first-order linear partial differential equation for $F$:  
\[
F_{x}= y F + z F_{y} + z F_{z} + w F_{w} + 
x w F_{x} - w y F_{y} - w z F_{z} - w^2 F_{w}.
\]
The solution of this PDE is a nice application of the standard method of characteristic curves, which yields 
\[
F(x,y,z,w)=e^{x(y-z)} \left( \frac{1 - w}{1 - we^{x(1 - w)}} \right)^{\frac{\textrm{\footnotesize\raisebox{1pt}{$z$}}}{\textrm{\footnotesize\raisebox{-1pt}{$w$}}}}.
\]
Putting $y=z=1$ in $F$  yields the generating function for free ascents:
\[
\left(\frac{1 - w}{1 - we^{(1 - w)x}} \right)^{\frac{\textrm{\footnotesize\raisebox{1pt}{$1$}}}{\textrm{\footnotesize\raisebox{-1pt}{$w$}}}}.
\]


\section{Further questions}

\vspace*{-3mm}

\hspace*{6mm} The fast recurrence for (1-23-4)-avoiding permutations established in this paper raises the question whether there might be fast but nonobvious recurrences for other similar patterns.  For example, there is a slow  recurrence for (12-34)-avoiding permutations posted to OEIS at 
\htmladdnormallink{A113226}{http://www.research.att.com:80/cgi-bin/access.cgi/as/njas/sequences/eisA.cgi?Anum=A113226}. Can you do better?

\end{document}